\documentclass[a4paper,12pt]{article}
\usepackage{amssymb}
\newtheorem{Th}{Theorem}

\newtheorem{Co}{Corollary}
\newtheorem{Lm}{Lemma}
\newtheorem{Lma}{Lemma}[section]

\newtheorem{Rm}{Remark}

\newcommand{\be}{\begin{equation}}
\newcommand{\ee}{\end{equation}}
\newcommand{\R}{\mathbb{R}}
\newcommand{\N}{\mathbb{N}}

\newcommand{\reset}{\setcounter{equation}{0}\setcounter{Th}{0}\setcounter{Prop}{0}\setcounter{Co}{0}
\setcounter{Lm}{0}\setcounter{Rm}{0}}

\def\lf{\left}
\def\rg{\right}

\def\la{\lambda}

\def\ep{\varepsilon}
\def\ds{\displaystyle}

\def\Om{\Omega}
\def\om{\omega}
\def\p{\partial}

\begin{document}
\title{Sub-criticality of Schr\"odinger Systems with Antisymmetric Potentials.}


\author{Tristan Rivi\`ere\footnote{Department of Mathematics, ETH Zentrum,
CH-8093 Z\"urich, Switzerland.}}
\date{ }
\maketitle

\noindent{\bf Abstract :} Let $m$ be an integer larger or equal to 3. We prove that  Schr\"odinger systems on $B^m$ with $L^{m/2}-$antisymmetric potential $\Omega$
of the form
\[
-\Delta v=\Omega\ v
\]
can be written in divergence form and we deduce that solutions $v$ in $L^{m/(m-2)}$ are in fact
$W^{2,q}_{loc}$ for any $q<m/2$.

\vspace{1cm}

\section{Introduction}

 In \cite{Ri1} the author proved the sub-criticality of the following linear systems in 2 Dimension
 \be
 \label{I.1}
 -\Delta u=\Omega\cdot \nabla u\quad,
 \ee
 where $u\in W^{1,2}(D^2,{\R}^n)$ and $\Om\in L^2(D^2,{\R}^2\otimes so(n))$ ($n$ is an arbitrary
 integer, $so(n)$ is the subspace of $M_n({\R})$, the  space of $n\times n$ square matrices, made of antisymmetric matrices) and we have using
 the matrix multiplication : in coordinates (\ref{I.1}) reads
 \[
 \forall i=1\cdots n\quad\quad-\Delta u^i=\sum_{j=1}^n\Omega_j^i\cdot\nabla u^j\quad.
 \]
 Precisely, it is proved in \cite{Ri1} that such a $u$ is in fact in $W^{2,p}_{loc}(D^2,{\R}^n)$ for
 every $p<2$. This result has been obtained by writing (\ref{I.1}) in conservative form. This
 was possible due to the following result
 \begin{Th}
 \label{th-I.1} \cite{Ri1}
 There exists a map, in a neighborhood of the origin, of the form 
 \be
 \label{I.2}
 \begin{array}{rcl}
\ds {\mathcal L}\  :\quad L^2(D^2,{\R}^2\otimes so(n))&\ds\longrightarrow & \ds L^\infty\cap W^{1,2}(D^2,Gl_n({\R}))\\[5mm]
\Om &\longrightarrow & A
\end{array}
\ee
such that
\be
\label{I.3}
div(\nabla_\Om A):=div(\nabla A-A\Om)=0\quad,
\ee
  and with the following control 
\be
\label{I.4}
\|dist(A,SO(n))\|_\infty+\|A\|_{W^{1,2}}\le C\|\Om\|_{L^2}\quad,
\ee
where $C$ is a positive constant independent of $\Om$.\hfill $\Box$
\end{Th}
Once $A$ is constructed one easily see that
\be
\label{I.5}
 -\Delta u=\Omega\cdot \nabla u\quad\quad\Longleftrightarrow\quad\quad  div(A\nabla u+B\nabla^\perp u)=0\quad .
 \ee
where $\nabla^\perp B:=(-\p_yB,\p_x B)=\nabla A-A\Om$. The higher integrability of $\nabla u$ is then a direct consequence of this conservative form 
of the system by applying Wente's estimates (see \cite{Ri1} and \cite{Ri2}).
This result has lead in particular to a proof of Hildebrandt's conjecture on the regularity
of critical points to conformally invariant problems in two dimension.

In this paper we will study this time Schr\"odinger systems of the form
\be
\label{I.6}
-\Delta v=\Om\  v\quad,
\ee
where $v\in L^{m/(m-2)}(B^m,{\R}^n)$ and $\Om\in L^{m/2}(B^m,so(n))$, $n$ is an arbitrary integer
and $m$ is  an arbitrary integer larger or equal to 3. $B^m_r$ denotes the $m-$dimensional ball centered
at the origin of ${\R}^m$ and when we don't write the subscript it implicitly means that $r=1$
(i.e $B^m$ denotes the unit ball).
In coordinates (\ref{I.6}) means
\[
 \forall i=1\cdots n\quad\quad-\Delta v^i=\sum_{j=1}^n\Omega_j^i\ v^j\quad.
 \]
Like (\ref{I.1}) in 2-dimension, the system (\ref{I.6}) is also a-priori critical
for $v\in L^{m/(m-2)}$ in $m$ dimension. Indeed, under these assumptions $v\in L^{m/(m-2)}$ and $\Om\in L^{m/2}$ we obtain that the r.h.s. of (\ref{I.6}) and hence
$\Delta v$ is in $L^1$ and, using classical singular integral theory, we deduce in return that
$v\in L^{m/(m-2),\infty}_{loc}$ which is ''almost'' the information we started from.
Such a structure in general situations offers no hope for having any of the properties that characterize
 sub-critical problems such as better
integrability of $v$, local uniqueness of the solutions...etc. It is {\it a-priori} simply critical. However, here again, the antisymmetry
of $\Omega$ will imply that sub-criticality in fact holds.

Our main result is the following.

\begin{Th}
\label{th-I.2}
Let $m\ge 3$ and $n\in{\N}^\ast$
 There exists a map, in a neighborhood of the origin, of the form 
 \be
 \label{I.7}
 \begin{array}{rcl}
\ds {\mathcal S}\  :\quad L^{m/2}(B^m, so(n))&\ds\longrightarrow & \ds L^\infty\cap W^{2,m/2}(B^m,Gl_n({\R}))\\[5mm]
\Om &\longrightarrow &A
\end{array}
\ee
such that
\be
\label{I.8}
\Delta A+A\Omega=0\quad .
\ee
and there exists $C>0$, independent of $\Omega$, such that
\be
\label{I.9}
\|dist(A,O(n))\|_{L^\infty(B^m_{1/2})}+\|A\|_{W^{2,m/2}(B^m)}\le C\|\Om\|_{L^{m/2}(B^m)}\quad.
\ee
\hfill $\Box$
\end{Th}
Once $A$ is constructed one easily observe that for any $v\in L^{m/(m-2)}$ the following equivalence holds
\be
\label{I.10}
-\Delta v=\Om\  v\quad\quad\Longleftrightarrow\quad\quad div(A\ \nabla v-\nabla A\ v)=0\quad.
\ee
We have then been able to write Schr\"odinger Systems with antisymmetric potential
in conservative form\footnote{Observe that the product $A\ \nabla v$ makes sense since $A\in W^{1,m}$, by Sobolev embeddings, and we have $A\ \nabla v:=\nabla(A\ v)-\nabla A\ v$.}.
A corollary of the existence of such conservation law for Schr\"odinger Systems with anti-symmetric potential is the sub-criticality of such systems. Precisely we have.
\begin{Co}
\label{co-I.3} Let $n\in{\N}^\ast$ and $m\ge 3$. Let $v\in L^{m/(m-2)}(B^m,{\R}^n)$ satisfying
\[
-\Delta v=\Om\  v\quad,
\]
where $\Om\in L^{m/2}(B^m,so(n))$, then $v\in W^{2,q}_{loc}(B^m,{\R}^n)$ for any $q<m/2$.
\hfill $\Box$
\end{Co}
Our results and their proofs take their source jointly in \cite{Ri1} but also in \cite{DR2} where
F. Da Lio and the author were studying the regularity of 1/2-harmonic maps from the real line into manifolds - see also \cite{DR1}.
They reduced the original problem to the one of proving that the following equation is sub-critical in one dimension
\[
\Delta^{1/4}v=\Omega\ v\quad,
\]
where $v\in L^2({\R},{\R}^n)$ and $\Om\in L^2({\R},so(n))$.
We end-up this introduction by making the following remarks.
\begin{Rm}
\label{rm-I.4}
It is important to insist on the fact that, a-priori, from the way we construct them, both the mappings ${\mathcal L}$ and ${\mathcal S}$ are not continuous between, respectively, 
 $L^2(D^2,{\R}^2\otimes so(n))$ and $L^\infty\cap W^{1,2}(D^2,Gl_n({\R})$ and between $L^{m/2}(B^m, so(n))$ and $L^\infty\cap W^{2,m/2}(B^m,Gl_n({\R}))$.
 Our constructions both in \cite{Ri1} and in the present paper are realized by the application of successively local inversion theorem and continuity argument like the construction of Coulomb Gauges for $L^{m/2}-$curvatures in \cite{Uh}. Recently a construction of ${\mathcal L}$
 using a more direct variational method has been proposed by A.Schikorra in \cite{Sc}. He was following 
 an approach introduced by F.H\'elein in order to construct ''Coulomb Moving Frames''
 (see \cite{He} lemma 4.1.3). A construction of ${\mathcal S}$ using such a variational argument
might a-priori be possible and would be interesting in itself. \hfill $\Box$
\end{Rm}
\begin{Rm}
\label{rm-I.5}
Though the two problems treated respectively in theorem~\ref{th-I.1} and theorem~\ref{th-I.2}
share many resemblances in the results, one of the main points which are given by  the $L^\infty-$control of $A$ in resp. (\ref{I.4}) and (\ref{I.9}) are obtained via two different arguments.
In the first problem the $L^\infty$-control of $A$ comes basically from the application of Wente estimates for Jacobian and the so-called ''integrability by compensation'' phenomenon, whereas in 
the second problem it comes from an application of the maximum principle.
This difference is very fundamental and striking at least to us.
\end{Rm}
\begin{Rm}
\label{rm-I.6}
In \cite{RS}, M. Struwe and the author established the sub-criticality of (\ref{I.1}) in arbitrary
dimension in Morrey spaces. This was motivated by applications to the partial regularity of
stationary critical points to conformally invariant Lagrangians in higher dimension. However the existence of the 
Matrix valued map $A$ in $L^\infty(B^m,Gl_n({\R}))$ satisfying
\[
div(\nabla_{\Omega}A)=0
\]
was problematic due to the fact that Wente integrability by compensation does not provide
$L^\infty$ bounds in the classical Morrey spaces but only in their Littlewood-Paley counterpart (see \cite{Ke}). Here however, since the $L^\infty$ control of $A$ in theorem~\ref{th-I.2} is obtained
by the application of the Maximum principle, the chances are high that theorem~\ref{th-I.2}
extends to higher dimension for the ad-hoc Morrey spaces which make system (\ref{I.6}) a-priori critical.
\end{Rm}
The paper is organized as follows. In section 2 we construct the map ${\mathcal S}$, proving
then theorem-\ref{th-I.2}, and using an intermediate construction of a solution $P\in W^{2,m/2}(B^m,SO(n))$ solving
\[
\frac{1}{2}\lf[\Delta P\ P^{-1}-P\ \Delta P^{-1}\rg]+P\ \Om\ P^{-1}=0
\]
that we postpone in the appendix. In section 3 we deduce from theorem~\ref{th-I.2} the corollary~\ref{co-I.3}.

\section{Proof of theorem~\ref{th-I.2}.}
\reset
Let $\Om\in L^{m/2}(B^m,so(n))$ and $v\in L^{m/{m-2}}(B^m,{\R}^n)$ satisfying (\ref{I.6}). Consider
$P\in W^{2,m/2}(B^m,SO(n))$ given by lemma~\ref{lm-A.1}.
We compute
\[
-\Delta(P\ v)=\Delta P\ v-P\ \Delta v-2\ div(\nabla P\ v)\quad.
\]
Introducing $w:=P\ v$, the equation (\ref{I.6}) is  then equivalent to
\[
-\Delta w=\lf[\Delta P\ P^{-1}+P\,\Omega\, P^{-1}\rg]\ w-2\ div(\nabla P\ P^{-1}\ w)\quad.
\]
Taking into account this special choice of $P$ we have made and satisfying (\ref{A.1}), with our notations
the system (\ref{I.6}) becomes equivalent to
\be
\label{II.1}
-\Delta w -\frac{1}{2}\lf[\Delta P\ P^{-1}+P\ \Delta P^{-1}\rg]\ w+2\ div(\nabla P\ P^{-1}\ w)=0\quad.
\ee
Observe that
\[
\begin{array}{rl}
\ds-\lf[\Delta P\ P^{-1}+P\ \Delta P^{-1}\rg]&\ds=-div(\nabla P\ P^{-1}+P\ \nabla P^{-1})+2\nabla P\cdot \nabla P^{-1}\\[5mm]
=-2(\nabla P\ P^{-1})^2
\end{array}
\]
where we have used twice that $\nabla P\ P^{-1}=-P\ \nabla P^{-1}$. The notation
for the r.h.s $-2(\nabla P\ P^{-1})^2$ has to be understood as follows
\[
-2(\nabla P\ P^{-1})^2:=-2\sum_{j=1}^m(\p_{x_j}P\ P^{-1})^2
\]
where the squares in the r.h.s refer to Matrix multiplication. Observe that each $\p_{x_j}P\ P^{-1}$
is an $L^m$ map taking values into $so(n)$ therefore each $-(\p_{x_j}P\ P^{-1})^2$
is an $L^{m/2}$ map taking values into the space $Sym^+_n({\R})$ of {\bf symmetric non-negative $n\times n-$matrices}\footnote{ Indeed if $a$ is a real antisymmetric matrix we have that $(a^2)^t=a^ta^t=a^2$ and for every x in ${\R}^n$ $<x,-(a)^2x>=-x^ta^2x=x^ta^tax=(ax)^tax\ge 0$}  . Hence
\[
-(\nabla P\ P^{-1})^2\in L^{m/2}(B^m,Sym_n^+({\R}))\quad.
\]
Combining (\ref{II.1}) with the previous observations, the Schr\"odinger system (\ref{I.6}) becomes equivalent to
\be
\label{II.2}
-\Delta w -(\nabla P\ P^{-1})^2\ w+2\ div(\nabla P\ P^{-1}\ w)=0\quad.
\ee
Standard elliptic estimates gives that for any given  $r<m/2$, if $\|\nabla P\|_{L^m}$ is small enough 
- depending on $r$ a-priori -, then there
exists a unique solution $Q\in W^{2,r}(B^m,M_n({\R}))$ of the following problem
\be
\label{II.3}
\lf\{
\begin{array}{l}
-\Delta Q-2\nabla Q\cdot\nabla P\ P^{-1}-Q\ (\nabla P\ P^{-1})^2=0\quad\quad\mbox{ in } B^m\\[5mm]
\quad Q=Id\quad\quad\mbox{ on }\p B^m
\end{array}
\rg.
\ee
This comes from the following a-priori estimates
\[
\|\nabla Q\cdot\nabla P\ P^{-1}\|_{L^r}\le\|\nabla Q\|_{L^{rm/m-r}}\ \|\nabla P\|_{L^m}\le C_r\ \|Q-Id\|_{W^{2,r}_0}\ \|\nabla P\|_{L^m}
\]
and 
\[
\begin{array}{rl}
\ds\|(Q-Id)\ (\nabla P\ P^{-1})^2\|_{L^r}&\ds\le  \|(Q-id)\|_{L^{rm/m-2r}}\    \|(\nabla P\ P^{-1})^2\|_{L^{m/2}}\\[5mm]
 &\le\ds C_r\ \|Q-Id\|_{W^{2,r}_0}\ \|\nabla P\|_{L^m}^2\quad .
 \end{array}
\]
We establish now the following lemma.
\begin{Lm}
\label{lm-II.1}
Let $m\ge 3$ and $n\in {\N}^\ast$. There exists $\ep_0>0$ such that for any $P\in W^{1,m}(B^m,SO(n))$ 
satisfying 
\[
\int_{B^m}|\nabla P|^m<\ep_0\quad,
\]
and any $Q\in W^{2,2m/(m+2)}(B^m,M_n({\R}))$ solving
\[
\lf\{
\begin{array}{l}
-\Delta Q-2\nabla Q\cdot\nabla P\ P^{-1}-Q\ (\nabla P\ P^{-1})^2=0\quad\quad\mbox{ in } B^m\\[5mm]
\quad Q=Id\quad\quad\mbox{ on }\p B^m\quad.
\end{array}
\rg.
\]
Then $Q\in L^\infty\cap W^{2,m/2}(B^m,M_n({\R}))$. Moreover there exists $C_m>0$ such that
\be
\label{II.4a}
\|dist(Q,O(n))\|_{L^\infty(B^m_{1/2})}\le C_m\ \lf[\int_{B^m}|\nabla P|^m\rg]^{2/m}\quad .
\ee
\hfill $\Box$
\end{Lm}
{\bf Proof of Lemma~\ref{lm-II.1}.}

We first show that for any $X\in{\R}^n$ the following inequality holds :
\be
\label{II.4}
\Delta(X^t\,Q\,Q^t\,X)\ge 0\quad .
\ee
We have
\[
\begin{array}{rl}
\ds\Delta(X^t\,Q\,Q^t\,X)&\ds=X^t\,\Delta Q\,Q^t\,X+X^t\,Q\,\Delta Q^t\, X+2X^t\,\nabla Q\cdot\nabla Q^t\, X\\[5mm]
\ds &\ds=-2\,X^t\,\nabla Q\cdot(\nabla P\ P^{-1})\,Q^t\,X-X^t\,Q(\nabla P\ P^{-1})^2\,Q^t\,X\\[5mm]
 &\ds\ +2\, X^t\, Q\ (\nabla P\ P^{-1})\cdot\nabla Q^t\, X-X^t\,Q\,(\nabla P\ P^{-1})^2\,Q^t\,X\\[5mm]
  &\ds +2X^t\,\nabla Q\cdot\nabla Q^t\, X
\end{array}
\]
where all this above operations make a distributional sense  (Leibnitz rule) as long as $Q\in W^{2,2m/(m+2)}(B^m)$, which is our assumption.
Observe that\footnote{Since for $Y$ and $Z$ in ${\R}^n$ we have $Y^t\,Z=Z^t\,Y$}
\[
\begin{array}{rl}
-2\,X^t\,\nabla Q\cdot(\nabla P\ P^{-1})\,Q^t\,X&\ds=-2\, ((\nabla P\ P^{-1})\,Q^t\,X)^t\cdot(X^t\,\nabla Q)^t\\[5mm]
 &\ds=2\, X^t\, Q\,(\nabla P\ P^{-1})\cdot\nabla Q^t\,X\quad.
 \end{array}
\]
Hence we have
\be
\label{II.5}
\begin{array}{rl}
\ds\Delta(X^t\,Q\,Q^t\,X)&\ds=+4\, X^t\, Q\ (\nabla P\ P^{-1})\cdot\nabla Q^t\, X\\[5mm]
 &\ -2X^t\,Q\,(\nabla P\ P^{-1})^2\,Q^t\,X +2X^t\,\nabla Q\cdot\nabla Q^t\, X\quad .
\end{array}
\ee
Cauchy-Schwartz inequality tells that
\[
\begin{array}{rl}
\ds-2\, X^t\, Q\ (\nabla P\ P^{-1})\cdot\nabla Q^t\, X&\ds\le X^t\, Q\ (\nabla P\ P^{-1})\cdot (\nabla P\ P^{-1})^t\, Q^t\, X\\[5mm]
 &\quad\ds+X^t\,\nabla Q\cdot\nabla Q^t\, X\quad .
 \end{array}
\]
Since again $(\nabla P\ P^{-1})^t=-(\nabla P\ P^{-1})$, the previous inequality implies
\be
\label{II.6}
\begin{array}{rl}
\ds 4\,X^t\, Q\ (\nabla P\ P^{-1})\cdot\nabla Q^t\, X &\ds\ge 2X^t\,Q\,(\nabla P\ P^{-1})^2\,Q^t\,X\\[5mm]
   &\ds\quad-2X^t\,\nabla Q\cdot\nabla Q^t\, X\quad.
\end{array}
\ee
Combining (\ref{II.5}) and (\ref{II.6}) we obtain (\ref{II.4}). Applying the Maximum Principle we obtain
that\footnote{Since $|Q^tX|^2=X^t\,Q\,Q^t\,X$.} 
\be
\label{II.7}
\sup_{X\in{\R}^n}\|Q^tX\|^2_{L^\infty(B^m)}\le 1\quad .
\ee
This implies that $Q\in L^{\infty}(B^m)$. Hence $Q\ (\nabla P\ P^{-1})^2\in L^{m/2}(B^m)$.
Since we have the {\it a-priori} estimate (for any $1<r<m$)
\[
\begin{array}{rl}
\|\nabla Q\cdot\nabla P\ P^{-1}\|_{L^r}&\ds \le \|\nabla P\|_{L^m}\ \|\nabla Q\|_{L^{rm/m-r}}\\[5mm]
 &\le C_r\epsilon_0\ \|Q-Id\|_{W^{2,r}_0(B^m)}\quad,
 \end{array}
 \]
Applying it successively for $r=2m/m+2$ and $r=m/2$ we deduce that, for $\epsilon_0$ chosen small
enough, the operator
\[
\begin{array}{rcl}
K_{P}\ :\ W^{2,r}_0(B^m,M_n({\R}))&\longrightarrow & L^r(B^m,M_n({\R}))\\[5mm]
 \eta&\longrightarrow &\ds -\Delta\eta-2\nabla\eta\cdot\nabla P\ P^{-1}
 \end{array}
 \]
is an isomorphism for both $r=2m/m+2$ and $r=m/2$. Applying it to $\eta=Q-Id$ we obtain, since 
$Q\ (\nabla P\ P^{-1})^2\in L^{m/2}(B^m)$, that $Q\in W^{2,{m/2}}(B^m, M_n({\R}))$ and the following
estimate holds
\be
\label{II.8}
\|Q-Id\|_{W^{2,m/2}_0(B^m)}\le C_m\ \lf[\int_{B^m}|\nabla P|^m\rg]^{2/m}\quad.
\ee
(\ref{II.4}) can also be written in the following way : $\forall X\in S^{n-1}$ - $S^{n-1}$ denotes the unit sphere.
\be
\label{II.9}
\lf\{
\begin{array}{l}
\ds\Delta(X^tX-X^t\,Q\,Q^t\,X)\le 0\quad\quad\mbox{ in }{\mathcal D}'(B^m)\\[5mm]
\ds X^tX-X^t\,Q\,Q^t\,X=0\quad\quad\mbox{ on }\quad\p B^m
\end{array}
\rg.
\ee
Hence we can apply Harnack Inequality to each function  $X^tX-X^t\,Q\,Q^t\,X\in L^{\infty}(B^m)$ for each $X\in S^{n-1}$ (see for instance \cite{GT}), and we have
\be
\label{II.10}
\begin{array}{rl}
\ds0\le\sup_{x\in B^m_{1/2}}X^tX-X^t\,Q\,Q^t\,X&\ds\le C_m\ \int_{B^m_{1/2}}X^tX-X^t\,Q\,Q^t\,X\\[5mm]
 &\ds = C_m\int_{B^m_{1/2}}(X^t-X^t\, Q)(X+Q\,X)\\[5mm]
 &\ds\le 2\, C_m\int_{B^m}|Q-Id|\\[5mm]
  &\ds\le C'_m \lf[\int_{B^m}|\nabla P|^m\rg]^{2/m}
\end{array}
\ee
where we used successively (\ref{II.7}) and (\ref{II.8}). Since we can exchange the sup quantificators, (\ref{II.10})  implies in particular
\be
\label{II.11}
\|\sup_{X\in S^{m-1}}|X^tX-|QX|^2|\|_{L^\infty( B^m_{1/2})}\le C'_m \lf[\int_{B^m}|\nabla P|^m\rg]^{2/m}\quad .
\ee
We have
 \[
 \begin{array}{rl}
\ds 2[X^tY-X^tQ^tQY]&\ds=(X+Y)^t(X+Y)-|Q(X+Y)|^2\\[5mm]
 &\ds +X^tX-|QX|^2+Y^tY-|QY|^2
 \end{array}
 \]
Hence 
\be
\label{II.12}
\|\sup_{X,Y\in S^{m-1}}|(X,Y)-(QX,QY)|\|_{L^\infty( B^m_{1/2})}\le C''_m \lf[\int_{B^m}|\nabla P|^m\rg]^{2/m}\quad .
\ee
where $(\cdot,\cdot)$ denotes the scalar product in ${\R}^n$.
Denote by $\||\cdot|\|$ the 2-norm on square matrices given by $\|| M\||^2=tr(M^tM)$. For $\epsilon_0$ sufficiently small (\ref{II.10}) implies that $Q$ is in a neighborhood of $O(n)$ in which the orthogonal projection $\pi_{O(n)}$ with respect to the scalar product $<M,N>:=tr(M^t N)$ onto $O(n)$ is smooth.
Denote $R:=\pi_{O(n)}(Q)$ and let $S:=R^{-1}(Q-R)$. $\||S\||=dist(Q,O(n))$. Because of the minimality property of $\||Q-R\||$ among all $R$ in $O(n)$, $S$ satisfies
\[
\forall a\in so(n), \quad\quad 0=<Ra,Q-R>=- tr(a\, R^{-1}(Q-R))=-tr(a\, S)
\]
which means that $S$ is symmetric $S^t=S$. Observe that
\be
\label{II.13}
\begin{array}{rl}
\ds\sup_{X,Y\in S^{m-1}}|(X,Y)-(QX,QY)|&=\sup_{X,Y\in S^{m-1}}|X^t(Id-Q^tQ)Y|\\[5mm]
 &\ds=\|Id-Q^tQ\|=\||2S+S^2\||
 \end{array}
\ee
Using the fact that $S$ is small in $L^\infty(B^m_{1/2})-$norm, we have $||2S+S^2\||\ge\||S\||$ and combining (\ref{II.12})
and (\ref{II.13}) we deduce (\ref{II.4a}) and lemma~\ref{lm-II.1} is proved. \hfill $\Box$

\medskip

\noindent{\bf End of the proof of Theorem~\ref{th-I.2}.}

We fix $2m/(m+2)=r$, we assume $\|\nabla P\|^m_{L^m}$ to be less than 
$\epsilon_0$ in lemma~\ref{lm-II.1} and we consider $Q$ given by this lemma.
Multiply (\ref{II.2}) on the left by $Q$ gives
\[
\begin{array}{l}
0=-Q\,\Delta w -Q\,(\nabla P\ P^{-1})^2\ w+2Q\ div(\nabla P\ P^{-1}\ w)\\[5mm]
\quad=-Q\,\Delta w-[Q\,(\nabla P\ P^{-1})^2+2\nabla Q\cdot\nabla P\ P^{-1}]\ w+2div(Q\ \nabla P\ P^{-1}\ w)\\[5mm]
\quad=-Q\Delta w+\Delta Q\, w+2div(Q\ \nabla P\ P^{-1}\ w)\\[5mm]
\quad=div(-Q\nabla w+\nabla Q\, w+2Q \nabla P\ P^{-1}\ w)
\end{array}
\]
Going back now to the original variable $v=P^{-1}w$ gives
\[
div((QP)\ \nabla v-\nabla (QP)\ v)=0
\]
and $A:=QP$ satisfies the conclusion of the theorem~\ref{th-I.2} which concludes the proof.\hfill $\Box$

\section{Proof of corollary~\ref{co-I.3}.}
\reset

Once we prove that $v$ belongs to $L^p_{loc}(B^m)$ for some $p>m/m-2$ a classical bootstrap
argument gives that $v\in W^{2,q}_{loc}(B^m)$ for any $q<m/2$.

In order to prove that $v$ belongs to $L^p_{loc}(B^m)$ for some $p>m/m-2$, it suffices to prove
that there exists $\gamma>0$ such that
\be
\label{III.1}
\sup_{x_0\in B^m_{1/2},r<1/4}r^{-\gamma}\lf[\int_{B^m_r(x_0)}|v|^{m/m-2}\rg]^{(m-2)/m}<+\infty\quad.
\ee
Indeed, this later fact injected in the system (\ref{I.6}) implies that
\be
\label{III.2}
\sup_{x_0\in B^m_{1/2},r<1/4}r^{-\gamma}\int_{B^m_r(x_0)}|\Delta v|<+\infty\quad,
\ee
and using Adams embedding results (see \cite{Ad}) one directly obtain that $v$ belongs to $L^p_{loc}(B^m_{1/2})$ for some $p>m/m-2$. Hence in order to prove corollary~\ref{co-I.3} it suffices to establish
a Morrey type estimate of the form (\ref{III.1}) that will be obtained by a very standard argument once
we use  theorem~\ref{th-I.2}.

Using theorem~\ref{th-I.2} we rewrite the system (\ref{I.6}) in the following form
\[
-\Delta w=-2div(\nabla A A^{-1} w)
\]
where $w=Av$. On the ball $B_r(x_0)$ we decompose $w=\phi+\xi$ where
\[
\lf\{
\begin{array}{l}
\ds-\Delta\phi=-2div(\nabla A A^{-1} w)\quad\quad\mbox{ in } B_r(x_0)\\[5mm]
\ds \phi=0\quad\quad\mbox{ on }\p B_r(x_0)
\end{array}
\rg.
\]
Hence $\xi$ is harmonic and for any $\la<1$ one has\footnote{Indeed, if $\xi$ is harmonic then
$\Delta|\xi|^{m/m-2}\ge 0$ and hence the function $\rho\rightarrow |B_\rho(x_0)|^{-1}\int_{B_{\rho}(x_0)}|\xi|^{m/m-2}$ is increasing.}
\be
\label{III.3}
\int_{B_{\la r}(x_0)}|\xi|^{m/m-2}\le \la^m \int_{B_r(x_0)} |\xi|^{m/m-2}\quad.
\ee
Standard elliptic estimate give for $\phi$
\be
\label{III.4}
\int_{B_r(x_0)}|\phi|^{m/m-2}\le C_m\|A^{-1}\|_\infty\ \lf[\int_{B_r(x_0)}|\nabla A|^m\rg]^{1/(m-2)}\int_{B_r(x_0)}|w|^{m/m-2}
\ee
where $C_m$ is independent of $r$. For any $\ep>0$ there exists a radius $r_\ep$ such that for any $x_0$ in $B_{1/2}^m$ and $r<r_\epsilon$ we have
$C_m\|A^{-1}\|_\infty\ \lf[\int_{B_r(x_0)}|\nabla A|^m\rg]^{1/(m-2)}<\ep$. We shall choose $\epsilon$ later.
Summing (\ref{III.3}) and (\ref{III.4}) gives
\be
\label{III.5}
\begin{array}{l}
2^{-2/(m-2)}\int_{B_{\la r}(x_0)}|w|^{m/m-2}\le \ \int_{B_{\la r}(x_0)}|\phi|^{m/m-2}+|\xi|^{m/m-2}\\[5mm]
\quad\quad\le \la^m \int_{B_r(x_0)} |w-\phi|^{m/m-2}+\ep\int_{B_r(x_0)}|w|^{m/m-2}\\[5mm]
\quad\quad\le (2^{2/m-2}\la^m+\ep+\la^m\ep 2^{2/m-2} )\int_{B_r(x_0)}|w|^{m/m-2}
\end{array}
\ee
We choose now $\la$ and $\ep$ small enough in such a way that
\[
2^{2/(m-2)}(2^{2/m-2}\la^m+\ep+\la^m\ep 2^{2/m-2} )\le1/2
\]
This gives
\be
\label{III.6}
\int_{B_{\la r}(x_0)}|w|^{m/m-2}\le1/2\ \int_{B_r(x_0)}|w|^{m/m-2},
\ee
from which we deduce a Morrey estimate of the form (\ref{III.1}) for $w$, which itself finally implies (\ref{III.1}) for $v$. Corollary~\ref{co-I.3} is then proved.\hfill $\Box$
\reset
\appendix
\section{Appendix}

The appendix is devoted to the proof of the following lemma.

\begin{Lma}
\label{lm-A.1}
Let $m\ge 3$ and $n\in{\N}^\ast$. There exists $\ep_0>0$ and $C>0$ such that,
for any $\Om\in L^{m/2}(B^m,so(n))$ there exists $P\in W^{2,m/2}(B^m,SO(n))$
satisfying
\be
\label{A.1}
\lf\{
\begin{array}{l}
\frac{1}{2}\lf[\Delta P\ P^{-1}-P\ \Delta P^{-1}\rg]+P\ \Om\ P^{-1}=0\quad\quad\mbox{ in }\quad{\mathcal D}'(B^m)\\[5mm]
P=Id_{SO(n)}\quad\quad\mbox{ on }\quad{\mathcal D}'(B^m)
\end{array}
\rg.
\ee
and
\be
\label{A.2}
\|P-Id\|_{W^{2,m/2}_0(B^m)}\le C\ \|\Om\|_{L^{m/2}}\quad .
\ee
\hfill $\Box$
\end{Lma}

{\bf Proof of lemma~\ref{lm-A.1}.} We follow a similar approach to the one introduced in the appendix of \cite{Ri1} which was itself inspired by the work of K.Uhlenbeck \cite{Uh}.
Let $q>m/2$ and $\varepsilon>0$. Consider
$$
{\cal{U}}_{\varepsilon}^q=\left\{\Omega\in L^q(B^m,so(n))\ :~\int_{\R}|\Omega|^{m/2} dx<\varepsilon\right\}\,.
$$
{\it  {\bf Claim:} There exist $\varepsilon_0>0$ and $C>0$ such that
\[
{\cal{V}}_{\epsilon_0,C}^q:=\lf\{
\begin{array}{c}
 \Omega\in {\cal{U}}_{\varepsilon_0}^q~:  \mbox{there exits $P$ satisfying (\ref{A.1}) and (\ref{A.2})}\\[5mm]
 \mbox{ and }P=exp(U)\mbox{  with $||U||_{W^{2,q}_0(B^m)}\le C \|\Om\|_{L^q(B^m)} $}
\end{array}
\rg\}
\]
is open and closed in ${\cal{U}}_{\varepsilon_0}^q$ for the $L^q-$norm and thus
${\cal{V}}_{\varepsilon_0,C}^q\equiv {\cal{U}}_{\varepsilon_0}^q$ (since  ${\cal{U}}_{\varepsilon_0}^q$ is clearly path connected)\,.
\par
}
This claim implies lemma~\ref{lm-A.1}. Indeed, for this $\ep_0$ we consider $\Om\in L^{m/2}(B^m,so(n))$ such that $\|\Om\|_{L^{m/2}}<\ep_0$. By convolutions one gets a sequence of
maps $\Om_k\in {\cal{U}}_{\varepsilon}^q$ converging strongly to $\Om$ in $L^2$.

 Let $P_k\in
W^{2,q}(B^m,SO(n))$ given by the claim and satisfying both (\ref{A.1}) and (\ref{A.2}) for $\Om_k$. We can extract a subsequence that weakly converges in $W^{2,m/2}(B^m,SO(n))$ to a limit $P$
in $W^{2,m/2}(B^m,M_n({\R}))$ .

By lower semicontinuity of the $W^{2,m/2}-$norm under weak convergence and by Rellich compactness embedding, we deduce that $P$ satisfies (\ref{A.2}) and that $P$ takes values into
the rotations $SO(n)$. Again by compactness embedding we have that $P_k$ converges strongly
to $P$ in every $L^q$ for $q<+\infty$ and since $\Delta P_k$ converges weakly to $\Delta P$ in $L^{m/2}$ we pass easily to the limit in the equation (\ref{A.1}) and lemma-\ref{lm-A.1} is proved.

\medskip

It then remains to prove the claim.

\medskip 

\noindent{\bf Step 1  :} {\it For any $\varepsilon_0>0$ and $C>0$ ${\cal{V}}_{\varepsilon_0}^q$
is closed in ${\cal{U}}_{\varepsilon_0}^q$ }. The proof of this step follows one by one the argument we just used to prove that the claim implies lemma~\ref{lm-A.1}.

\medskip

It then remains to establish the following.

\medskip

\noindent{\bf Step 2  :} {\it There exists $\varepsilon_0>0$ and $C>0$ such that ${\cal{V}}_{\varepsilon_0,C}^q$
is open in ${\cal{U}}_{\varepsilon_0}^q$ }. 

Before to establish the step 2, we will prove a lemma that roughly tells us that as soon
as $\|P-Id\|_{W^{2,m/2}}$ is small enough then (\ref{A.2}) automatically holds. Precisely we have.

\begin{Lma}
\label{lm-A.2}
Let $m\ge 3$ and $n\in {\N}^\ast$. There exists $\ep_1>0$ and $C_1>0$ such that for any
$P\in W^{2,m/2}(B^m, SO(m))$ sucht that $P=Id$ on $\p B^m$, if
\be
\label{A.2a}
\|P-Id\|_{W^{2,m/2}_0(B^m)}\le \ep_1
\ee
then
\be
\label{A.3}
\|P-Id\|_{W^{2,m/2}_0(B^m)}\le C_1\ \| P^{-1}\ \Delta P-\Delta P^{-1}\ P\|_{L^{m/2}(B^m)}\quad ,
\ee
and such that for any $P\in W^{2,q}(B^m, SO(m))$ satisfying $P=Id$ on $\p B^m$ and (\ref{A.2a})
we have also
\be
\label{A.3a}
\|P-Id\|_{W^{2,q}_0(B^m)}\le C_1\ \| P^{-1}\ \Delta P-\Delta P^{-1}\ P\|_{L^{q}(B^m)}\quad .
\ee
\hfill$\Box$
\end{Lma}
{\bf Proof of lemma~\ref{lm-A.2}.}
We write
\be
\label{A.4}
\begin{array}{rl}
\ds P^{-1}\Delta P&\ds=\frac{1}{2}\lf[P^{-1}\Delta P-\Delta P^{-1}\ P\rg]\\[5mm]
 &\ds\ +\frac{1}{2}\lf[P^{-1}\ \Delta P+\Delta P^{-1}\ P\rg]
 \end{array}
\ee
Moreover we have 
\be
\label{A.5}
\begin{array}{rl}
\ds P^{-1}\ \Delta P+\Delta P^{-1}\ P&\ds= div\lf(P^{-1}\nabla P+\nabla P^{-1}\ P\rg)-2\nabla\, P^{-1}\cdot\nabla P\\[5mm]
  &\ds=-2\nabla\, P^{-1}\cdot\nabla P
  \end{array}
 \ee 
Hence, by assumption, we have
\be
\label{A.6}
\begin{array}{rl}
\|P^{-1}\ \Delta P+\Delta P^{-1}\ P\|_{L^{m/2}(B^m)}&\le 2\|\nabla P\|_{L^m(B^m)}\ \|\nabla P\|_{L^m(B^m)}\\[5mm]
 &\ds \le 2\ep_1\ \|\nabla P\|_{L^m(B^m)}
 \end{array}
 \ee
Since $P-Id=0$ on $\p B^m$, standard elliptic estimates give
\[
\|\nabla P\|_{L^m(B^m)}\le C_m\ \|\Delta P\|_{L^{m/2}(B^m)}\quad.
\]
This last fact combined with (\ref{A.5}) and (\ref{A.6}) give for $2\ep_1\ C_m<1/2$
\[
\|\Delta P\|_{L^{m/2}(B^m)}\le \frac{2}{3}\|P^{-1}\Delta P-\Delta P^{-1}\ P\|_{L^{m/2}(B^m)}\quad .
\]
Using again the fact that  $P-Id=0$ on $\p B^m$, standard elliptic estimates combined with the previous inequality gives (\ref{A.3}).

(\ref{A.3a}) is proved in a similar way. Observe that

\be
\label{A.6z}
\begin{array}{rl}
\|P^{-1}\ \Delta P+\Delta P^{-1}\ P\|_{L^{q}(B^m)}&\le 2\|\nabla P\|_{L^m(B^m)}\ \|\nabla P\|_{L^{qm/m-q}(B^m)}\\[5mm]
 &\ds \le 2\ep_1\ \|\nabla P\|_{L^{qm/m-q}(B^m)}
 \end{array}
 \ee
Since $P-Id=0$ on $\p B^m$, standard elliptic estimates give
\[
\|\nabla P\|_{L^{qm/m-q}(B^m)}\le C_m\ \|\Delta P\|_{L^{q}(B^m)}\quad.
\]
and we finish the argument as in the case $q=m/2$ in order to get (\ref{A.3a}) this completes the proof of lemma~\ref{lm-A.2}. \hfill $\Box$

We start now the proof of step 2. Intrduce the map  $F$ defined as follows
\[
\begin{array}{rcl}
\ds F\ :\ W^{2,q}_0(B^m,so(n))&\longrightarrow &\ds L^q(B^m,so(n))\\[5mm]
U\ &\longrightarrow & P^{-1}\Delta P-\Delta P^{-1}\ P
\end{array}
\]
where $P=exp\,(U)$. We first prove that the map $F$ is $C^1$. This comes from the following facts
\begin{itemize}
\item[i)] Since $W^{2,q}$ for $q>m/2$ embedds continuously in $C^0$, the map $U\rightarrow exp\,(U)$ is clearly smooth from $W^{2,q}_0(B^m,so(n))$ into $W^{2,q}(B^m,SO(n))$.
\item[ii)] The operator $\Delta$ is a smooth linear map from $W^{2,q}(B^m,M_n({\R}))$ into
$L^q(B^m, M_n({\R}))$.
\item[iii)] Since again $W^{2,q}$ embedds continuously in $L^\infty$ - $W^{2,q}$ is an algebra -
the following map
\[
\begin{array}{rcl}
\ds \Pi\ :\ W^{2,q}_0(B^m,M_n({\R}))\times L^q(B^m,M_n({\R}))&\longrightarrow & L^q(B^m,M_n({\R}))\\[5mm]
(A,B)\quad&\longrightarrow&\quad A\,B
\end{array}
\]
is also smooth.

\end{itemize}
 For $v$ and $w$ in $so(n)$,  we denote
 $$
 D(v)\cdot w:=exp(-v)\ d\, exp_{v}\cdot w:=exp(-v)\ \frac{d}{dt}exp(v+tw)_{|_{t=0}}\in so(n)\quad .
 $$
With this notation we have
\[
\begin{array}{l}
\ds dF_{U_0}\cdot\eta=(P_0\ D(U_0)\cdot\eta)^t\ \Delta P_0-\Delta P_0^{-1}\ (P_0\ D(U_0)\cdot\eta)\\[5mm]\ds\quad\quad+P_0^{-1}\ \Delta(P_0\ D(U_0)\cdot\eta)-\Delta(P_0\ D(U_0)\cdot\eta)^t\ P_0
\end{array}
\]
where $P_0:=exp(U_0)$. Denote $\zeta:=D(U_0)\cdot\eta$ for $\eta\in W^{2,q}_0(B^m,so(n))$, we observe that
\be
\label{A.6a}
\begin{array}{l}
\ds \frac{1}{2}dF_{U_0}\cdot\eta=L_{P_0}\zeta:=\Delta\zeta+[P_0^{-1}\nabla P_0,\nabla\zeta]+[\Omega_0,\zeta]
\end{array}
\ee
where $2\Omega_0:=P_0^{-1}\ \Delta P_0-\Delta P_0^{-1}\ P_0$
We now establish the following lemma

\begin{Lma}
\label{lm-A.3}
There exists $\varepsilon_2>0$ such that for any $U_0\in W^{2,q}_0(B^m,so(n))$ satisfying
\be
\label{A.7}
\|exp(U_0)-Id\|_{W^{2,m/2}}\le\varepsilon_2\quad,
\ee
then $dF_{U_0}$ is invertible between $W^{2,q}_0(B^m,so(n))$ and $L^q(B^m,so(n))$.
\hfill$\Box$
\end{Lma}
{\bf Proof of Lemma~\ref{lm-A.3}.}
We first prove that there exists $\varepsilon>0$ such that whenever $\|exp(U_0)-Id\|_{W^{2,m/2}}\le\varepsilon$, there exists $C_{U_0}>0$, such that for any $\om\in L^q(B^m,so(n))$ there exists
a unique $\zeta\in W^{2,q}_0(B^m,so(n)) $ for which
\be
\label{A.8}
\lf\{
\begin{array}{l}
\ds L_{P_0}\zeta=\om\quad,\\[5mm]
\ds\|\zeta\|_{W^{2,q}_0(B^m,so(n))}\le C_0\ \|\om\|_{ L^q(B^m,so(n))}
\end{array}
\rg.
\ee
Since $W^{2,q}_0(B^m)$ embedds continuously in $L^\infty(B^m)$ it is clear that
$[\Omega_0,\zeta]\in L^q$. Moreover 
\[
\|[P_0^{-1}\nabla P_0,\nabla\zeta]\|_{L^q}\le 2\|\nabla P_0\|_{L^m}\ \|\nabla\zeta\|_{L^{qm/m-q}}
\le C\|P_0-id\|_{W^{2,m/2}_0}\ \|\zeta\|_{W^{2,q}_0}\quad .
\]
Hence $L_{P_0}$ is sending continuously $W^{2,q}_0(B^m,so(n))$ into $L^q(B^m,so(n))$.
Since $m>q>m/2$ we have that $4/m-1/q>2/m$. We can hence choose $r$ such that $4/m-1/q>1/r>2/m$
(for instance ${1}/{r}:=3/m-1/2q$). For such a $r$ we have
\be
\label{A.8a}
\|[\Om_0,\zeta]\|_{L^r}\le 2\ \|\Om_0\|_{L^{m/2}}\ \|\zeta\|_{L^{rm/m-2r}}\le C_q\ \|\Om_0\|_{L^{m/2}}\ \|\zeta\|_{W^{2,r}_0}
\ee
and 
\be
\label{A.8b}
\|[P_0^{-1}\nabla P_0,\nabla\zeta]\|_{L^r}\le 2\|\nabla P_0\|_{L^m}\ \|\nabla\zeta\|_{L^{rm/m-r}}
\le C\|P_0-id\|_{W^{2,m/2}_0}\ \|\zeta\|_{W^{2,r}_0}\quad .
\ee
Hence using standard elliptic theory, we obtain that for $\|P_0-id\|_{W^{2,m/2}_0}$ small enough, for
any $\om\in L^r(B^m,M_n({\R}))$
there exists a unique solution $\zeta$ in $W^{2,r}_0(B^m,M_n({\R})$ of $L_{P_0}\zeta=\om$.
Assume moreover that $\om$ takes values into $so(n)$ then we have, since $(P_0^{-1}\ \nabla P_0)^t=-P_0^{-1}\ \nabla P_0$ and $\Om_0^t=-\Om_0$,
\[
L_{P_0}(\zeta+\zeta^t)=0\quad.
\]
The just proved uniqueness result gives then $\zeta^t=-\zeta$. 

Hence we have established that
\[
\begin{array}{rcl}
L_{P_0}\ :\ W^{2,r}_0(B^m,so(n)) &\longrightarrow & L^r(B^m, so(n))\\[5mm]
\zeta\ &\longrightarrow & \Delta\zeta+[P_0^{-1}\nabla P_0,\nabla\zeta]+[\Omega_0,\zeta]
\end{array}
\]
is an isomorphism.

Let $1/s:=1/q+1/r-2/m$. Our assumption on $r$ gives $1/s<2/m$. Denoting $\Delta^{-1}_0$ the Inverse
of the laplacian on $B^m$ for the zero Dirichlet boundary data, we have
\[
\|\Delta^{-1}_0([\Om_0,\zeta])\|_\infty\le C\ \|[\Om_0,\zeta]\|_{L^s}\le C\ \|\Om_0\|_{L^q}\ \|\zeta\|_{L^{mr/m-2r}}\le C\ \|\Om_0\|_{L^q}\ \|\zeta\|_{W^{2,r}_0}\quad.
\] 
Moreover
\[
\begin{array}{rl}
\ds\|\Delta^{-1}_0([P_0^{-1}\nabla P_0,\nabla\zeta])\|_\infty&\ds\le C\ \|[P_0^{-1}\nabla P_0,\nabla\zeta]\|_{L^s}\le C\
\|\nabla P_0\|_{L^{qm/m-q}}\ \|\nabla \zeta\|_{L^{rm/m-r}}\\[5mm]
 &\le C\ \|P_0-Id\|_{W^{2,q}_0}\ \|\zeta\|_{W^{2,r}_0}
 \end{array}
\]
From the two previous estimates we deduce that for any $\om\in L^q(B^m,so(n))$, the unique solution
$\zeta\in W^{2,r}_0(B^m,so(n))$ of $L_{P_0}\zeta=\om$ is in fact in $L^\infty$ and the following estimate
holds
\be
\label{A.9}
\begin{array}{rl}
\ds\|\zeta\|_{L^\infty(B^m)}&\ds\le C_q\ \|P_0-Id\|_{W^{2,q}_0(B^m)}\ \|\zeta\|_{W^{2,r}_0(B^m)}+C_q\ \|\om\|_{L^q(B^m)}\\[5mm]
 &\ds\le C_q\ \lf[1+\|P_0-Id\|_{W^{2,q}_0(B^m)}\rg]\ \|\om\|_{L^q(B^m)}
 \end{array}
 \ee
We then obtain that 
\be
\label{A.10}
\|[\Om_0,\zeta]\|_{L^q(B^m)}\le C_q\ \|\Delta P_0\|_{L^q}\ \lf[1+\|P_0-Id\|_{W^{2,q}_0(B^m)}\rg]\ \|\om\|_{L^q(B^m)}
\ee
Observe that inequality (\ref{A.8b}) is valid for any $r<m$ and hence in particular it holds for $q$ :
we have for any $\xi$ in $W^{2,q}_0$
\[
|[P_0^{-1}\nabla P_0,\nabla\xi]\|_{L^q}\le 2\|\nabla P_0\|_{L^m}\ \|\nabla\zeta\|_{L^{qm/m-q}}
\le C\|P_0-id\|_{W^{2,m/2}_0}\ \|\zeta\|_{W^{2,q}_0}\quad .
\] 
Hence for $\|P_0-id\|_{W^{2,m/2}_0}$ having been chosen small enough, by standard
elliptic estimates, the following map
\[
\begin{array}{rcl}
H_{P_0}\ :\ W^{2,q}_0(B^m,so(n)) &\longrightarrow & L^q(B^m, so(n))\\[5mm]
\xi\ &\longrightarrow & \Delta\xi+[P_0^{-1}\nabla P_0,\nabla\xi]
\end{array}
\]
is an isomorphism. Let $\xi:=H_{P_0}^{-1}\lf[\om-[\Om,\zeta]\rg]$. 
The argumentation we followed above for $L_{P_0}$ applies to $H_{P_0}$ in order
to show that it realizes an isomorphism between $W^{2,r}_0(B^m,so(n))$ and $L^r(B^m,so(n))$.
Hence since $H_{P_0}(\xi-\zeta)=0$ we deduce that $\zeta=\xi$ and hence we have proved that
$\zeta\in W^{2,q}_0(B^m,so(n))$ and the following estimate holds :
\[
\|\zeta\|_{W^{2,q}_0(B^m,so(n))}\le C_q\ \lf[1+\|\Delta P_0\|_{L^q}\ \lf[1+\|P_0-Id\|_{W^{2,q}_0(B^m)}\rg]\rg]\  \|\om\|_{L^q(B^m)}
\]
We have then established (\ref{A.8}). The map
\[
\begin{array}{rcl}
D\ :\ so(n)&\longrightarrow& Gl(so(n))\subset so(n)\otimes(so(n))^\ast\\[5mm]
v&\longrightarrow& exp(-v)\ dexp_v
\end{array}
\]
is smooth, hence, since $W^{2,q}(B^m)$ is an algebra for $q>m/2$ and since it embedds
in $C^0$ the map $D(U_0)$ realizes an isomorphism from $W^{2,q}_0(B^m, so(n))$ into
itself. Since now $dF_{U_0}=2\,L_{P_0}\circ D(U_0)$, we have proved lemma~\ref{lm-A.3}.
\hfill $\Box$.

\medskip

{\bf End of the proof of step 2.} We fix an $\varepsilon_0$ smaller than the $\varepsilon_1$
 of lemma~\ref{lm-A.2} and smaller than the $\varepsilon_2$ of lemma~\ref{lm-A.3}. Consider also
 $C$ equal to $C_1$ given by lemma~\ref{lm-A.2}. Let $\Om_0\in {\mathcal V}_{\varepsilon_0,C}^q$.
 According to Lemma~\ref{lm-A.3} we can apply the local inversion theorem and then there exists
 a neighborhood of $\Om_0$ in $L^q(B^m,so(n))$ such that for any $\Omega$ in this neighborhood
 there exists $P\in W^{2,q}(B^m,SO(n))$ such that (\ref{A.1}) holds. In particular this is true for any $\Omega$ in the intersection of this neighborhood with ${\mathcal U}_{\varepsilon_0}^q$. Since $\varepsilon_0\le\varepsilon_1$, Lemma~\ref{lm-A.2} applies and we deduce that all these $\Omega$
 belong to ${\mathcal V}_{\varepsilon_0,C}$. Hence we have proved that there exists a neighborhood
 of $\Omega_0$ whose intersection with ${\mathcal U}^q_{\varepsilon_0}$ is included in ${\mathcal V}_{\varepsilon_0,C}^q$. This shows that for this choice of $\varepsilon_0$ and $C$ ${\mathcal V}_{\varepsilon_0,C}^q$ is open in ${\mathcal U}_{\varepsilon_0}^q$. We have the proved step 2
 and we deduce lemma~\ref{lm-A.1}.\hfill $\Box$


\begin{thebibliography}{99}
 \bibitem[Ad]{Ad} Adams, David R. ''A note on Riesz potentials.'' Duke Math. J. 42 (1975), no. 4, 765--778.
\bibitem[DR1]{DR1}  Da Lio Francesca, Rivi\`ere Tristan ''3-commutator estimates and the regularity of 1/2-harmonic maps into spheres'', preprint arXiv:0901.2533, (2009).
\bibitem[DR2]{DR2}  Da Lio Francesca, Rivi\`ere Tristan ''The regularity of solutions to critical non-local Schr\"odinger systems on the line with antisymmetric potential and applications'' in preparation (2009).
\bibitem[GT]{GT} Gilbarg, David, Trudinger, Neil S. ''Elliptic Partial Differential Equations of Second Order''
Grundlehren der mathematischen Wissenschaften, 224, Springer (2001).
\bibitem[He]{He} H\'elein, Fr\'ed\'eric ``Harmonic Maps, Conservation Laws, and Moving Frames." Cambridge Tracts in Mathematics, 150. Cambridge University Press (2002).
\bibitem[Ke]{Ke} Keller, Laura ''Integrability by compensation in Besov-Morrey spaces and applications''
preprint (2009).
\bibitem[Ri1]{Ri1} Rivi\`ere, Tristan ``Conservation laws for conformally invariant variational problems." Invent. Math., 168 (2006), no 1, 1-22.
\bibitem[Ri2]{Ri2} Rivi\`ere, Tristan ''Integrability by Compensation in the Analysis of Conformally Invariant Problems.'' Minicourse at Pisa, to appear in Edizioni della Normale.
\bibitem[RS]{RS} Rivi\`ere, Tristan, Struwe, Michael ''Partial regularity for harmonic maps and related problems.'' Comm. Pure Appl. Math. 61 (2008), no. 4, 451--463.
\bibitem[Sc]{Sc} Schikorra, Armin ''A Remark on Gauge Transformations and the
Moving Frame Method'' preprint arXiv:0906.1972 (2009).
\bibitem[Uh]{Uh} Uhlenbeck Karen. ''Connections with $L^p$ bounds on curvature.''
Comm. Math. Phys., 83:31--42, 1982.



\end{thebibliography}
\end{document}